\def\ann{\mathop{\rm ann}\nolimits}
\def\ara{\mathop{\rm ara}\nolimits}
\def\Ass{\mathop{\rm Ass}\nolimits}

\def\Hom{\mathop{\rm Hom}\nolimits}
\def\InjH{\mathop{\rm E}\nolimits}
\def\LCMo{\mathop{\rm H}\nolimits}

\def\min{\mathop{\rm min}\nolimits}
\def\Spec{\mathop{\rm Spec}\nolimits}

\def\Naturalsign{{\rm l\kern-.23em N}}
\font \normal=cmr10 scaled \magstep0 \font \mittel=cmr10 scaled \magstep1 \font \gross=cmr10 scaled \magstep5
\input amssym.def
\input amssym.tex
\parindent=0pt
\gross Local Cohomology and complete intersections of arithmetic rank one
\normal
\bigskip
Michael Hellus
\par
hellus@math.uni-leipzig.de
\bigskip
\mittel
0. Introduction
\normal
\bigskip
In algebraic geomtry there is the important concept of complete intersections;
roughly speaking a complete intersection is a variety cut out by codimension
many equations; in general the number of equations needed to cut out a variety
given by an ideal is called the arithmetic rank (of the ideal). It was shown
in [8] that the notion of arithmetic rank is strongly related to the concept of
regular sequences on the Matlis duals of certain local cohomology modules
(basics on regular sequences can be found in [5], [11] and [12]) and thus is
related to the set of associated primes of such Matlis duals.
\par
While in [8]
the general situation with no restrictions on arithmetic rank and
cohomological dimension was
investigated, this work concentrates on the case ``cohomological dimension of
the given ideal is (at most) one'' and also on the interplay between
non-graded local and graded situations. The main results are new
characterizations for the case $\ara \leq 1$ (theorem 1 and theorem 2), the fact
that given a graded ring $R$ and a homogenous ideal $I$ of cohomological
dimension at most one the inclusion
$$\{ x\in I\vert x\hbox { non-homogenous}\} \subseteq \bigcup _{\goth p\in
\Ass _R(D(\LCMo ^1_I(R)))}\goth p$$
holds (theorem 4, $D$ is a Matlis dual functor) and a (somewhat technical) statement on regular sequences on
Matlis duals of certain local cohomology modules (theorem 5), which shows that
regular sequences on such Matlis duals behave well in some sense although
these modules are not finite in general.
\par
Section 1 describes the local (non-graded) situation, section 2 the graded
situation, while section 3 contains the results and section 4 presents some
open questions.
\bigskip
\mittel
1. The local situation
\normal
\bigskip
Given a noetherian local ring $(R,\goth m)$ and an ideal $I$ of $R$, the
arithmetic rank of $I$ is defined as the minimal number of generators of $I$
up to radical:
$$\ara (I):=\min \{ l\in \Naturalsign \vert \exists r_1,\dots r_l\in R\vert
\sqrt {I}=\sqrt {(r_1,\dots ,r_l)R}\} \ \ . $$
If we denote by $\LCMo ^s _I(R)$ ($s\in \Naturalsign $) the local cohomology
modules of $R$ supported in $I$ then apparantely (e. g. by Koszul cohomology
arguments) the implication $\ara (I)\leq 1 \Rightarrow 0=\LCMo ^2_I(R)=\LCMo
^3_I(R)=\dots $ holds, while the converse implication ``$\Leftarrow $'' does not
hold in general (counterexample: $k$ a field, $S=k[[x_1,x_2,x_3,x_4]]$ a
formal power series ring over $k$ in variables $x_1,x_2,x_3,x_4$, $\goth
p:=\sqrt {(f:=x_1x_4-x_2x_3,x_2^4-x_1^3x_4,x_3^4-x_4^3x_1)S}, R:=S/fS,
I:=\goth p/fS$. It is easy to see that both
$$0=\LCMo ^2_I(R)=\LCMo ^3_I(R)=\dots $$
and $\ara (I)\geq 2$ hold. The variety defined by $\goth p$ is a
variant of the so-called ``Macaulay-curve''). We have the equivalence
$$\ara (I)\leq 1\iff 0=\LCMo ^2_I(R)=\LCMo ^3_I(R)=\dots \hbox { and }\exists
f\in I:f\hbox { operates surjectively on } \LCMo ^1_I(R)\ \ .$$
The (easy) proof can be extracted from [8, section 0], but we
repeat it here for sake of completeness:
\par
$\Rightarrow $: Assume $\sqrt I=\sqrt {fR}$ for some $f\in R$. W. l. o. g. we can assume $f\in I$. $f$ operates surjectively on
$\LCMo ^1_{fR}(R)=\LCMo ^1_I(R)$. $\Leftarrow $: $\LCMo ^1_I(\_)$ is a
right-exact functor (on $R$-modules). Therefore we have an exact sequence
$$\LCMo ^1_I(R)\buildrel f\over \to \LCMo ^1_I(R)\to \LCMo ^1_I(R/fR)\to 0\ \
.$$
Thus $\LCMo ^1_I(R/fR)=0$ holds, implying $\LCMo ^1_I(R/\goth p)=0$ for all
$\goth p\in \Spec (R)$ containing $f$. $\sqrt I=\sqrt {fR}$ follows.
\par
From now on we will always assume $0=\LCMo ^2_I(R)=\LCMo ^3_I(R)=\dots $
unless stated otherwise. We denote by $\InjH :=\InjH _R(R/\goth m)$ the
$R$-injective hull of $R/\goth m$ and by $D$ the Matlis dual functor from
$(R-mod)$ to itself sending $M$ to $\Hom _R(M,\InjH )$. Then we have
$$\ara (I)\leq 1 \iff \exists f\in I: f\hbox { operates injectively on }
D(\LCMo ^1_IR())\iff I\not \subseteq \bigcup _{\goth p\in \Ass _R(D(\LCMo
^1_I(R)))}\goth p\ \ .$$
\bigskip
\mittel
2. The graded situation
\normal
\bigskip
Next we consider the following situation (referred to from now on as ``graded
situtation''): $K$ a field, $l\in \Naturalsign ,R=K[X_0,\dots ,X_N]/J$ ($N\in
\Naturalsign ,J\subseteq K[X_0,\dots ,X_n]$ a homogenous ideal, where every
$X_i$ has a multidegree in $\Naturalsign ^l$), $I\subseteq R$ a homogenous
ideal, $\goth m$ the maximal homogenous ideal $(X_0,\dots ,X_N)R$ of $R$,
$\InjH :=\InjH _R(R/\goth m)$ an $R$-injective hull of $R/\goth m$; $\InjH $
has a natural
grading and serves also as a *-$R$-injective hull of $R/\goth m$ (for details
on *-notation see [3, sections 12 and 13]). $*D$ shall denote the functor from graded $R$-modules
to itself sending $M$ to $*\Hom _R(M,\InjH )$. The homogenous arithmetic
rank of $I$ is defined as
$$\ara ^h:=\min \{ l\in \Naturalsign \vert \exists r_1,\dots ,r_l\in R\hbox {
homogenous}:\sqrt I=\sqrt {(r_1,\dots ,r_l)R}\} \ \ .$$
Just like in the non-graded local case one can show
$$\eqalign {\ara ^h(I)\leq 1&\iff \exists f\in I\hbox { homogenous}: f\hbox { operates
injectively on }(*D)(\LCMo ^1_I(R))\cr &\iff I^h\not \subseteq \bigcup _{\goth p\in
\Ass _R((*D)(\LCMo ^1_I(R)))}\goth p\ \ .\cr }$$
Here $I^h:=\{ r\in I\vert r\hbox { homogenous}\} $ and all $\goth p\in
\Ass _R((*D)(\LCMo ^1_I(R)))$ are homogenous.
\bigskip
\mittel
3. Results
\normal
\bigskip
Let $(R,\goth m)$ be a noetherian local ring and $X\subseteq \Spec (R)$ a
subset. We say ``$X$ satisfies prime avoidance'' if, for every ideal
$J$ of $R$,
$$J\subseteq \bigcup _{\goth p\in X}\goth p$$
implies
$$\exists \goth p_0\in X:J\subseteq \goth p_0\ \ .$$
{\bf Theorem 1}
\smallskip
Let $(R,\goth m)$ be a noetherian local ring and $I$ an ideal of $R$ such that
$0\buildrel (*)\over =\LCMo ^2_I(R)=\LCMo ^3_I(R)=\dots $. Then
$$\ara (I)\leq 1\iff \Ass _R(D(\LCMo ^1_I(R)))\hbox { satisfies prime
avoidance}\ \ .$$
Proof:
\par
We set
$$D:=D(\LCMo ^1_I(R))\ \ .$$
\par
$\Rightarrow $: Let $J\subseteq R$ be an ideal such that
$$J\subseteq \bigcup _{\goth p\in \Ass _R(D)}\goth p\ \ .$$
Assuming $\Hom _R(R/J,D)=0$ we conclude $\LCMo ^1_I(R/J)=0$ (by Matlis
duality). Because of $(*)$ it follows that $\LCMo ^1_I(R/\goth p)=0$ for all
prime ideals $\goth p$ of $R$ containing $J$. Again because of $(*)$ we have
$I\subseteq \goth p$ for all $\goth p$ containing $J$, that is $I\subseteq \sqrt
J$. There is an $x\in R$ such that $\sqrt I=\sqrt {xR}$. Hence $x^l\in J$ for
$l>>0$. So there is a $\goth p\in \Ass _R(D)$ containing $x$. Now we have
$$0=\LCMo ^1_{xR}(R/\goth p)=\LCMo ^1_I(R/\goth p)$$
and thus
$$0=D(\LCMo ^1_I(R/\goth p))=\Hom _R(R/\goth p,D(\LCMo ^1_I(R)))$$
contradicting $\goth p\in \Ass _R(D)$. The assumption $\Hom _R(R/J,D)=0$ is false and
so there exists $d\in D\setminus \{ 0\} $ such that $J\subseteq \ann _R(d)$.
\par
$\Leftarrow $:
We have to show the existence of an $x\in I$ operating surjectively on
$\LCMo^1_I(R)$. Assume to the contrary
$$I\subseteq \bigcup _{\goth p\in \Ass_R(D)}\goth p\ \ .$$
From the hypothesis we get a $\goth p_0\in \Ass _R(D)$ such that $I\subseteq
\goth p_0$. But this $\goth p_0$ would satisfy
$$0\neq \LCMo ^1_I(R/\goth p_0)=0\ \ .$$
\smallskip
Similarly, in the graded situation, let $X\subseteq \Spec ^h(R):=\{ \goth p\in
\Spec (R)\vert \goth p\hbox { homogenous}\} $ be any subset. $X$ satisfies
``homogenous prime avoidance'' if, for every homogenous ideal $J$ of $R$,
$$J^h\subseteq \bigcup _{\goth p\in X}\goth p$$
implies
$$\exists \goth p_0\in X: J\subseteq \goth p_0\ \ .$$
\bigskip
{\bf Theorem 2}
\smallskip
Let $R$ be graded and $I\subseteq R$ an homogenous ideal such that
$0=\LCMo ^2_I(R)=\LCMo ^3_I(R)=\dots $. Then
$$\ara ^h(I)\leq 1\iff \Ass _R((*D)(\LCMo ^1_I(R)))\hbox { satisfies homgenous
prime avoidance}$$
holds.
\par
Proof:
\par
The proof consists mainly of a graded version of the proof of theorem 1:
\par
$\Rightarrow $: Let $J\subseteq R$ be an homogenous ideal such that
$J^h\subseteq \bigcup _{\goth p\in \Ass _R((*D)(\LCMo ^1_I(R)))}\goth p$ and
$x\in R^h$ an element such that $\sqrt I=\sqrt {xR}$. We assume
$$\Hom _R(R/J,*\Hom _R(\LCMo ^1_I(R),\InjH ))=0$$
and remark that for the first $\Hom $ (in the preceeding formula) it would not
make any difference if we replaced $\Hom $ by $*\Hom $. This implies
$$*\Hom _R((R/J)\otimes _R\LCMo ^1_I(R),\InjH )=0$$
and hence $\LCMo ^1_I(R/J)=0$. Thus $I\subseteq \goth q$ for all prime ideals $\goth
q$ of $R$ containing $J$. This implies the existence of a $\goth p_0\in \Ass
_R((*D)(\LCMo ^1_I(R)))$ such that $x\in \goth p_0$ contradicting $\LCMo
^1_I(R/\goth p_0)\neq 0$.
\par
$\Leftarrow $: We assume that for every $x\in I^h$ there exists a $\goth p\in
\Ass _R((*D)(\LCMo ^1_I(R)))$ such that $x\in \goth p$, i. e.
$$I^h\subseteq \bigcup _{\goth p\in \Ass _R(*\Hom _R(\LCMo ^1_I(R),\InjH
))}\goth p\ \ .$$
There is a $\goth p_0\in \Ass _R(*\Hom _R(\LCMo ^1_I(R),\InjH ))$ containing
$I$, contradicting $\LCMo ^1_I(R/\goth p_0)\neq 0$.
\bigskip
{\bf Remark}
\smallskip
In the graded situation, given graded $R$-modules $M$ and $N$,
$$*\Hom _R(M,N)\subseteq \Hom _R(M,N)$$ holds. For finite $M$ one has
equality here, but for arbitrary $M$ equality does not hold in general. In fact one has
$$\Ass _R(*\Hom _R(M,N))\subsetneq \Ass _R(\Hom _R(M,N))$$ in general as we
will see below (in the remark following theorem 4) in the case $M=\LCMo
^1_I(R)$, $N=\InjH $; then we will also see that (in
some sense) $\Ass _R(\Hom _R(\LCMo ^1_I(R),\InjH))$ is much larger than $\Ass _R(*\Hom
_R(\LCMo ^1_I(R),\InjH))$.
\par
Still in the graded situation let $I$ be an ideal of $R$ (such that $0=\LCMo
^2_I(R)=\LCMo ^3_I(R)=\dots $). For every $f\in I$ we have
$$\sqrt I=\sqrt {fR}\iff \forall _{\goth p\in \Ass _R(\Hom _R(\LCMo ^1_I(R),\InjH ))} f\not \in
\goth p\ \ .$$
We replace the last condition by
$$\forall _{\goth p\in \Ass _R(*\Hom _R(\LCMo ^1_I(R),\InjH ))} f\not \in
\goth p$$
and get a weaker condition, which we will denote by $^h\kern-.12em\sqrt
I=^h\kern-.52em\sqrt {fR}$. Our next result shows there are many (inhomogenous) $f\in I$ such that
$^h\kern-.12em\sqrt  I=^h\kern-.47em\sqrt {fR}$.
\bigskip
{\bf Theorem 3}
\smallskip
Let $I$ be a homogenous ideal of $R$ such that $\ara ^h(I)\leq 1$. Let
$g_1,\dots ,g_n\in I\setminus \{ 0\} $ be homogenous of pairwise different degrees (in
$\Naturalsign ^l$) and such that
$$\sqrt I=\sqrt {(g_1,\dots ,g_n)R}\ \ .$$
Then
$$^h\kern-.12em\sqrt I=^h\kern-.46em\sqrt {(g_1+\dots +g_n)R}$$
holds.
\par
Proof:
\par
We have $\LCMo ^1_I(R/(g_1,\dots ,g_n)R)=0$ and hence
$$(g_1,\dots ,g_n)R\not\subseteq \goth p$$
for all $\goth p\in \Ass _R(*\Hom _R(\LCMo ^1_I(R),\InjH ))$. Theorem 2
implies
$$((g_1,\dots ,g_n)R)^h\not \subseteq \bigcup _{\goth p\in \Ass _R(*\Hom
_R(\LCMo ^1_I(R),\InjH ))} \goth p\ \ .$$
Because of the different degrees of the $g_i$ we conclude
$$(g_1+\dots +g_n)R\not \subseteq \bigcup _{\goth p\in \Ass _R(*\Hom
_R(\LCMo ^1_I(R),\InjH ))} \goth p$$
and the statement follows.
\bigskip
{\bf Lemma}
\smallskip
Let $R$ be a graded domain and $f\in R\setminus \{ 0\} $. Then the ideal
$\sqrt {fR}$ is homogenous if and only if $f$ is homogenous. In particular,
for any homogenous ideal $I$ fo $R$ we have
$$(\min \{ l\in \Naturalsign \vert \exists r_1,\dots r_l\in R\vert
\sqrt {I}=\sqrt {(r_1,\dots ,r_l)R}\}=:)\ara (I)\leq 1\iff \ara ^h(I)\leq 1\ \
.$$
\par
Proof:
\par
$\Leftarrow $ is clear. $\Rightarrow $: $R$ is $\Naturalsign ^l$-graded. This
given grading may be seen as $l$ given $\Naturalsign $-gradings on $R$ and so
we may assume $l=1$. Let $\delta :=\deg (f)$. Then $f_\delta (=$degree-$\delta
$-part of $f) \in \sqrt {fR}$, i. e. $\exists n\in \Naturalsign ^+$ and
$\exists g\in R: f_\delta ^n=fg$. $R$ is a domain and so $f$ (as well as $g$)
must be homogenous.
\smallskip
The preceeding lemma implies that if $R$ is a graded domain and $I\subseteq R$
is a homogenous ideal such that $\ara (I)\leq 1$ ($\iff \ara ^h(I)\leq 1$),
every non-homogenous $f\in I$ does not operate injectively on $\Hom _R(\LCMo
^1_I(R),\InjH )$. Furthermore, if $\ara (I)>1$ ($\iff \ara ^h(I)>1$), it is
clear by the remarks preceeding theorem 2, that no $f\in I$ operates
injectively on $\Hom _R(\LCMo ^1_I(R),\InjH )$. Thus we have the following
(somewhat surprising) result:
\bigskip
{\bf Theorem 4}
$$(\{ x\in I\vert x\hbox { non-homogenous} \} =:)I^{nh}\subseteq \bigcup _{\goth p\in \Ass _R(\Hom _R(\LCMo
^1_I(R),\InjH ))}\goth p$$
{\bf Remark}
\smallskip
While theorem 3 says there are (many) non-homogenous $f\in I$ operating
injectively on $*\Hom _R(\LCMo
^1_I(R),\InjH )$, theorem 4 says there are no non-homogenous $f\in I$
operating injectively on $\Hom _R(\LCMo
^1_I(R),\InjH )$.
\bigskip
Now we consider a more general situation: Let $(R,\goth m)$ be a noetherian
local ring, $I\subseteq R$ an ideal, $h\in \Naturalsign $ and we assume $\LCMo
^l_I(R)\neq 0\iff l=h$ holds. For every $R$-regular sequence $r_1,\dots
,r_h\in I$ the equivalence
$$\sqrt I=\sqrt {(r_1,\dots ,r_h)R}\iff r_1,\dots ,r_h \hbox { is a regular
sequence on } D(\LCMo ^h_I(R))$$
holds ([8, section 0]). But as $D(\LCMo ^h_I(R))$ is not finitely generated, regular sequences
on $D(\LCMo ^h_I(R))$ do not behave very well; e. g. one cannot expect that
all maximal regular sequences on $D(\LCMo ^h_I(R))$ have the same length. But now
we will see that at least some properties which are valid for regular
sequences in the finite
case remain true for $D(\LCMo ^h_I(R))$.
\par
The following fact is well-known: If $M$ is a finite $R$-module and $r_1,\dots ,r_h\in R$ is a $M$-regular
sequence then $r_1^\prime ,\dots ,r_h^\prime \in R$ is also an $M$-regular
sequence provided $(r_1^\prime ,\dots ,r_h^\prime )R=(r_1,\dots ,r_h)R$. In
our case it is clear that if a $R$-regular sequence $r_1,\dots ,r_h\in I$ is a
$D(\LCMo ^h_I(R))$-regular sequence then a $R$-regular-sequence $r_1^\prime
,\dots ,r_h^\prime \in I$ is also $D(\LCMo ^h_I(R))$-regular if $(r_1^\prime
,\dots ,r_h^\prime
)R=(r_1,\dots ,r_h)R$ holds (simply because of $\sqrt {(r_1^\prime ,\dots
,r_h^\prime )R}=\sqrt I$). But more is true:
\bigskip
{\bf Theorem 5}
\smallskip
Let $1\leq h^\prime \leq h$ and let $r_1,\dots ,r_{h^\prime }\in I$ be a
$R$-regular sequence that is $D(\LCMo ^h_I(R))$-regular. Furthermore, let
$r_1^\prime ,\dots ,r_{h^\prime }^\prime \in I$ be such
that $(r_1^\prime , \dots ,r_{h^\prime }^\prime )R=(r_1,\dots ,r_{h^\prime
})R$ holds. Then $r_1^\prime ,\dots ,r_{h^\prime }^\prime $ is a $D(\LCMo ^h_I(R))$-regular
sequence. In particular, any permutation of $r_1,\dots ,r_{h^\prime }$ is
again a $D(\LCMo ^h_I(R))$-regular sequence.
\par
Proof:
\par
It is clear that $r_1^\prime ,\dots ,r_{h^\prime }^\prime $ is an $R$-regular
sequence. By induction on $s\in \{ 1,\dots ,h^\prime \} $ we show two statements:
$$\LCMo ^l_I(R/(r_1,\dots ,r_s)R)\neq 0\iff l=h-s$$
and
$$D(\LCMo ^{h-s}_I(R/(r_1,\dots ,r_s)R))=D(\LCMo ^h_I(R))/(r_1,\dots
,r_s)D(\LCMo ^h_I(R))\ \ :$$
$s=1$: The short exact sequence
$$0\to R\buildrel r_1\over \to R\to R/r_1R\to 0$$
induces a short exact sequence
$$0\to \LCMo ^{h-1}_I(R/r_1R)\to \LCMo ^h_I(R)\buildrel r_1\over \to \LCMo
^h_I(R)\to 0$$
and shows
$$\LCMo ^l_I(R/r_1R)\neq 0\iff l=h-1\ \ .$$
The statement $D(\LCMo ^{h-1}_I(R/r_1R))=D(\LCMo ^h_I(R))/r_1D(\LCMo ^h_I(R))$
follows easily.
\par
$s>1$: The short exact sequence
$$0\to R/(r_1,\dots ,r_{s-1})R\buildrel r_s\over
\to R/(r_1,\dots ,r_{s-1})R\to R/(r_1,\dots ,r_s)R\to 0$$
induces an exact
sequence
$$0\to \LCMo ^{h-s}_I(R/(r_1,\dots ,r_s)R)\to \LCMo ^{h-(s-1)}_I(R/(r_1,\dots
,r_{s-1})R)\buildrel r_s\over \to \LCMo ^{h-(s-1)}_I(R/(r_1,\dots ,r_{s-1})R)$$
By induction hypothesis, $D(\LCMo ^{h-(s-1)}_I(R/(r_1,\dots
,r_{s-1})R))=D(\LCMo ^h_I(R))/(r_1,\dots ,r_{s-1})D(\LCMo ^h_I(R))$ and so, by
assumption, $r_s$ operates surjectively on $\LCMo ^{h-(s-1)}_I(R/(r_1,\dots
,r_{s-1})R)$ and we get
$$\LCMo ^l_I(R/(r_1,\dots ,r_s)R)\neq 0\iff l=r-s$$
and
$$\eqalign {D(\LCMo ^{h-(s-1)}_I(R/(r_1,\dots
,r_{s-1})R))&=D(\LCMo ^{h-(s-1)}_I(R/(r_1,\dots ,r_{s-1})R))/r_sD(\LCMo
^{h-(s-1)}_I(R/(r_1,\dots ,r_{s-1})R))\cr &=D(\LCMo ^h_I(R))/(r_1,\dots
,r_s)D(\LCMo ^h_I(R))\ \ .\cr }$$
In particular for $s=h^\prime $ we have
$$\LCMo ^l_I(R/(r_1,\dots ,r_{h^\prime })R)\neq 0\iff l=h-h^\prime \ \ ,$$
i. e.
$$\LCMo ^l_I(R/(r_1^\prime ,\dots ,r_{h^\prime }^\prime )R)\neq 0\iff
l=h-h^\prime \ \ .$$
We prove by descending induction on $s\in \{ 0,\dots ,h^\prime -1\} $ three
statements:
$$r_{s+1}^\prime \hbox { operates surjectively on } \LCMo
^{h-s}_I(R/(r_1^\prime ,\dots ,r_s^\prime )R)\ \ ,$$
$$\LCMo ^{h-l}_I(R/(r_1^\prime ,\dots ,r_s^\prime )R)\neq 0\iff l=s$$
and
$$D(\LCMo ^{h-(s+1)}_I(R/(r_1^\prime ,\dots ,r_{s+1}^\prime )R))=D(\LCMo
^{h-s}_I(R/(r_1^\prime ,\dots ,r_s^\prime )R))/r_{s+1}^\prime D(\LCMo
^{h-s}_I(R/(r_1^\prime ,\dots ,r_s^\prime )R))\ \ :$$
$s=h^\prime -1$: We consider the long exact $\Gamma _I$-sequence belonging to
the short exact sequence
$$0\to R/(r_1^\prime ,\dots ,r_{h^\prime -1}^\prime )R\buildrel r_{h^\prime
}^\prime \over \to R/(r_1^\prime ,\dots ,r_{h^\prime -1}^\prime )R\to
R/(r_1^\prime ,\dots ,r_{h^\prime }^\prime )R\to 0\ \ :$$
The surjectivity of $r_{h^\prime }^\prime $ on $\LCMo ^{h-(h^\prime
-1)}_I(R/(r_1^\prime ,\dots ,r_{h^\prime -1}^\prime )R)$ is obvious and the
other statements follow from the fact that for $l\neq h-(h^\prime -1)$ we have
injectivity of $r_{h^\prime }^\prime $ on $\LCMo ^l_I(R/(r_1^\prime ,\dots
,r_{h^\prime -1}^\prime )R)$, hence
$$\LCMo ^l_I(R/(r_1^\prime ,\dots
,r_{h^\prime -1}^\prime )R)=0$$
as $r_{h^\prime }^\prime \in I$.
\par
$s<h^\prime -1$: We consider the long exact $\Gamma _I$-sequence belonging to
the short exact sequence
$$0\to R/(r_1^\prime ,\dots ,r_s^\prime )R\buildrel r_{s+1}^\prime \over \to
R/(r_1^\prime ,\dots ,r_s^\prime )R\to R/(r_1^\prime ,\dots ,r_{s+1}^\prime
)R\to 0\ \ .$$
Now one can use similar arguments like in the case $s=h^\prime -1$ to prove
all three statements.
\bigskip
\mittel
4. Open questions
\normal
\bigskip
- Do there exist versions of theorem 1 and theorem 2 for the case $\ara (I)\leq l$ for
arbitrary $l$?
\smallskip
- In the remark preceeding theorem 3 the relation $^h\kern-.12em\sqrt I={^h\kern-.12em\sqrt {fR}}$ was
  defined (in a special situation). Is there a natural generalization of this
  to define relations like $^h\kern-.12em\sqrt I=^h\kern-.49em\sqrt {(f_1,\dots ,f_l)R}$?
\smallskip
- Suppose in the situation of theorem 1 that $X:=\Ass _R(D(\LCMo ^1_I(R)))$
  does not satisfy prime avoidance. The proof of theorem 1 shows that $I$ is contained in
  $\bigcup _{\goth p\in Ass_R(D(\LCMo ^1_I(R)))}\goth p$ and so there are
  ideals $K$ containing $I$ and being itself contained in $\bigcup _{\goth p\in
  Ass_R(D(\LCMo ^1_I(R)))}\goth p$ and being maximal with these properties. How
  are these ideals $K$ related to $I$?
\def\litem{\par\noindent \hangindent=\parindent\ltextindent}
\def\ltextindent#1{\hbox to \hangindent{#1\hss}\ignorespaces}
\bigskip
\mittel
\vfil
\eject
{\bf References}
\normal
\smallskip
\parindent=0.8cm
\litem{1.} Bass, H. On the ubiquity of Gorenstein rings, {\it Math. Z.} {\bf
82}, (1963) 8--28.
\medskip
\litem{2.} Brodmann, M. and Hellus, M. Cohomological patterns of coherent
sheaves over projective schemes, {\it Journal of Pure and Applied Algebra}
{\bf 172}, (2002) 165--182.
\medskip
\litem{3.} Brodmann, M. P. and Sharp, R. J. Local Cohomology, {\it Cambridge
studies in advanced mathematics} {\bf 60}, (1998).
\medskip
\litem{4.} Bruns, W. and Herzog, J. Cohen-Macaulay Rings, {\it Cambridge
University Press}, (1993).
\medskip
\litem{5.} Eisenbud, D. Commutative Algebra with A View Toward Algebraic
Geometry, {\it Springer Verlag}, (1995).
\medskip
\litem{6.} Grothendieck, A. Local Cohomology, {\it Lecture Notes in
Mathematics, Springer Verlag}, (1967).
\medskip
\litem{7.} Hellus, M. On the set of associated primes of a local cohomology
module, {\it J. Algebra} {\bf 237}, (2001) 406--419.
\medskip
\litem{8.} Hellus, M. On the associated primes of Matlis duals of top local cohomology
modules, to appear in {\it Communications in Algebra} {\bf 33}.
\medskip
\litem{9.} Huneke, C. Problems on Local Cohomology, {\it Res.
Notes Math. } {\bf 2}, (1992) 93--108.
\medskip
\litem{10.} Matlis, E. Injective modules over Noetherian rings, {\it Pacific J. Math.} {\bf 8}, (1958) 511--528.
\medskip
\litem{11.} Matsumura, H. Commutative ring theory, {\it Cambridge
University Press}, (1986).
\medskip
\litem{12.} Scheja, G. and Storch, U. Regular Sequences and Resultants, {\it AK
Peters}, (2001).

\end